\documentclass[conference]{IEEEtran}
\IEEEoverridecommandlockouts
\usepackage{cite}
\usepackage{amsmath,amssymb,amsfonts,amsthm}
\usepackage[flushleft]{threeparttable}

\usepackage{algorithmic}
\usepackage{graphicx}
\usepackage{textcomp}
\usepackage{xcolor}
\def\BibTeX{{\rm B\kern-.05em{\sc i\kern-.025em b}\kern-.08em
    T\kern-.1667em\lower.7ex\hbox{E}\kern-.125emX}}

\usepackage{footmisc}
\usepackage{url}

\usepackage{tabularray}
\usepackage{booktabs}

\usepackage{bm}

\usepackage{comment}

\makeatletter
\newcommand{\newlineauthors}{%
  \end{@IEEEauthorhalign}\hfill\mbox{}\par
  \mbox{}\hfill\begin{@IEEEauthorhalign}
}
\makeatother

\begin{document}

\title{Impact Analysis of Optimal EV Bi-directional Charging with Spatial-temporal Constraints\\
\thanks{This research is a collaborative effort with Enzen and received funding from the RACE for 2030 Cooperative Research Centre.}
}

\author{\IEEEauthorblockN{Xian-Long Lee}
\IEEEauthorblockA{\textit{Faculty of Information Technology} \\
\textit{Monash University}\\
Melbourne, Australia \\
xian-long.lee@monash.edu}
\and
\IEEEauthorblockN{Adel N. Toosi}
\IEEEauthorblockA{\textit{Faculty of Engineering and IT} \\
\textit{The University of Melbourne}\\
Melbourne, Australia \\
adel.toosi@unimelb.edu.au}
\and
\IEEEauthorblockN{Peter Pudney}
\IEEEauthorblockA{\textit{Industrial AI Research Centre} \\
\textit{University of South Australia}\\
Adelaide, Australia \\
peter.pudney@unisa.edu.au}
\newlineauthors
\IEEEauthorblockN{Ian McLeod}
\IEEEauthorblockA{\textit{Smart Energy} \\
\textit{Enzen Australia}\\
Adelaide, Australia \\
ian.mcleod@enzen.com}
\and
\IEEEauthorblockN{Muhammad Aamir Cheema}
\IEEEauthorblockA{\textit{Faculty of Information Technology} \\
\textit{Monash University}\\
Melbourne, Australia \\
aamir.cheema@monash.edu}
\and
\IEEEauthorblockN{Hao Wang}
\IEEEauthorblockA{\textit{Faculty of Information Technology} \\
\textit{Monash University}\\
Melbourne, Australia \\
hao.wang2@monash.edu}
}

\maketitle

\begin{abstract}
The growth in Electric Vehicle (EV) market share is expected to increase power demand on distribution networks. Uncoordinated residential EV charging, based on driving routines, creates peak demand at various zone substations depending on location and time. Leveraging smart charge scheduling and Vehicle-to-Grid (V2G) technologies offers opportunities to adjust charge schedules, allowing for load shifting and grid support, which can reduce both charging costs and grid stress. In this work, we develop a charge scheduling optimization method that can be used to assess the impact of spatial power capacity constraints and real-time price profiles. We formulate a mixed-integer linear programming problem to minimize overall charging costs, taking into account factors such as time-varying EV locations, EV charging requirements, and local power demands across different zones. Our analysis uses real data for pricing signals and local power demands, combined with simulated data for EV driving plans. Four metrics are introduced to assess impacts from the perspectives of both EV users and zones. Results indicate that overall EV charging costs are only minimally affected under extreme power capacity constraints.
\end{abstract}



\section{Introduction}
The trend of rising Electric Vehicle (EV) adoption and uncoordinated charging patterns is expected to add significant power demand to grid infrastructure \cite{iea2024globalEV, muratori2018impact}. EVs traveling throughout a metropolitan area create a clustering effect in local power demand, impacting various zones within the distribution network, e.g. a substation. If unconstrained, these charging demands from EVs can add significant stress to the grid. 
In metropolitan areas, EVs will be parked for significant durations at locations such as at home, at transit stations, near workplaces, or at shopping centres, according to users’ schedules and routines. While parked at various locations throughout the day, EVs can perform destination charging \cite{yong2023electric} in other locations in addition to charging at their residence. These charging activities occur at different times and locations throughout the day, increasing charging demand across zones during different periods.

Different charging periods are expected to have different prices due to Time-of-Use (ToU) pricing. Additionally, power constraints in different local demand zones can limit the power available for EV charging. Researchers have investigated the impact of EV charging in the distribution network under ToU pricing \cite{jones2022impact}. Another study analyzed the effects of implementing ToU pricing compared to uncoordinated EV charging within the distribution network \cite{sohail2023impact}. Additionally, with the growing adoption of EVs and their charging requirements, studies demonstrate the impact on power demand and increased stress on the grid \cite{powell2022charging}. Various studies evaluate different aspects of the distribution network \cite{flataker2022impact, rahman2020review, rahman2022comprehensive}. Results indicate an increase in charging demand in work areas and a reduction in residential areas. These findings have motivated us to conduct further impact analysis with power constraints and ToU price profiles to assess their effects on EV charging strategies.


In this paper, we develop an optimization method that can be used to assess the impact of pricing profiles and zone power constraints on uni-directional and bi-directional (or Vehicle-to-grid, V2G) charge scheduling optimization problems. We demonstrate how optimal charging strategies are influenced by ToU pricing and zone power constraints. The model minimizes overall costs while accounting for various factors, such as pricing profiles, EV driving preferences, and local power demand and constraints in different zones. Specifically, the key contribution of this research work aims to conduct in-depth analysis of three key aspects: 1) uni-directional and bi-directional charging, 2) pricing profiles, and 3) zone power constraints. To evaluate the analysis, we introduce four metrics to evaluate the impact performance from two perspectives: the zone's perspective and the EV user's perspective. These four metrics are calculated for the two analysis schemes in the optimization model. 

\section{System model}

\subsection{Framework}

\begin{figure}[ht]
\centering
\includegraphics[width=\columnwidth]{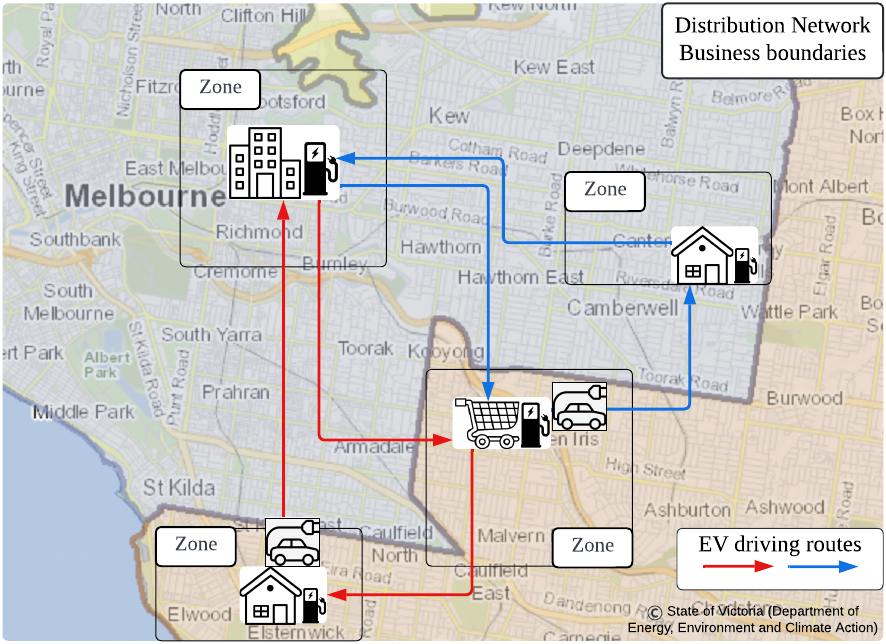}
\caption{Illustration of the scope of our work}\label{network}
\end{figure}

Figure \ref{network} illustrates EVs traveling throughout a metropolitan area based on their daily driving routines. Residential EVs that drive to different locations, such as workplaces, supermarkets, and homes, are placed across multiple ``zones'' within the metropolitan area. A ``zone'' represents an abstract area, which can be defined by geographic location or spatial grouping based on specific interests.
These zones may contain a mix of homes, workplaces, and shopping centers where EVs are parked for varying durations, creating potential charging needs. As EVs frequently travel across multiple zones, we expect an increase in charging demand that clusters with each zone’s local power demand. Additionally, zones may have distinct power constraints. To address this, the model minimizes overall charging costs, determining when and how much energy to charge while meeting each zone’s power constraints.


In our research, we aim to minimize the total charging costs for a group of EVs owned by individuals. The formulation seeks to reduce charging expenses within the EVs' charge schedules while incorporating zone-specific power constraints. 
We adopt a centralized approach to calculate the overall increase in charging demand from this group of EVs. Our experimental analysis focuses on the impacts of zone power capacity constraints and EV user perspectives. We evaluate the total power demand aggregated from both the EVs' charging requirements and local demand across different zones. 

\subsection{Problem Formulation}
In this section we formulate the objective function and constraints. The optimization is over a set $\mathbb{T}$ of equal time intervals, a set $\mathbb{Z}$ of zones, and a set $\mathbb{I}$ of EVs. Each time interval has duration $\Delta t$. 
The objective is to minimize the overall cost for the fleet of EVs, as defined in \eqref{obj}. EV users specify their availability for charging and discharging using the binary parameters $b_{i,z,t} \in \{0, 1\}$, where $b_{i,z,t} = 1$ signifies that EV $i$ will be connected in zone $z$ at time $t$. The variable $p_{i,t}$ is the controlled average charge or discharge power for EV $i \in \mathbb{I}$ during time interval $t \in \mathbb{T}$, and $\lambda(z, t, p)$ is the price of electricity in zone $z$ during time interval $t$, which depends on whether the EV is charging ($p > 0$) or discharging $(p < 0)$ as defined in constraint~\eqref{con_price}, where $\lambda_{z, t}$ is the charge price in zone $z$ during time interval $t$ and $\eta_{z, t}$ is the discharge price. The controlled average power is restricted by the charging point's maximum charging ($p^{+}$) and discharging ($p^{-}$) in constraint~\eqref{con_power}.
When EV $i$ is not be connected then constraint~\eqref{con_nopow} ensures that $p_{i,t} = 0$.
\begin{align}
    \min_{p_{i,t}} & \sum_{i\in\mathbb{I}} \sum_{z \in \mathbb{Z}} \sum_{t\in\mathbb{T}} b_{i,z,t} \, \lambda(z, t, p_{i, t})\, p_{i,t}\, \Delta{t} \label{obj} \\
    & \text{subject to} \notag \\
    & {p}^{-} \leq {p}_{i,t} \leq {p}^{+}, \forall t \in \mathbb{T}, i \in \mathbb{I} \label{con_power} \\
    & \lambda(z, t, p) =
      \begin{cases}
          \lambda_{z,t}, & \text{if}\ p \ge 0 \\
          \eta_{z,t}, & \text{if}\ p < 0
      \end{cases} \label{con_price} \\
    & p_{i,t} = 0 \quad \text{if} \quad \sum_{z \in \mathbb{Z}} b_{i,z,t} = 0 \label{con_nopow} \\
    & e_{i,t} = e_{i,t-1} + p_{i,t}\Delta{t} - d_{i,t} \label{con_bat} \\
    & 0 \leq {e}_{i,t} \leq e^\text{cap}_{i} \label{con_cap} \\
    & l_{z,t} + \sum_{i\in\mathbb{I}} p_{i,t} \, b_{i,z,t} \leq c^{+}_{z,t}, z \in \mathbb{Z}, t \in \mathbb{T} \label{con_consume} \\
    & e^\text{ini}_{i} + \sum_{t\in\mathbb{T}} \left(p_{i,t}\Delta{t} - v_{i,t}\right) \geq e^\text{tgt}_{i}, i \in \mathbb{I}. \label{con_tgt}
\end{align}

The amount of energy $e_{i,t}$ in EV $i$ at the end of time interval~$t$ is calculated by constraint~\eqref{con_bat}, where $d_{i,t}$ is the energy consumed while driving. Note that we assume the average power is charged or discharged directly to or from the battery, and the efficiency factor is not within the scope of this research. Each EV's battery energy is also subject to the battery capacity constraint~\eqref{con_cap} where $e^\text{cap}_i$ is the battery capacity of EV $i$. Constraint~\eqref{con_consume} ensures that in each zone $z$ and for every time interval $t$, the sum of the local energy demand $l_{z,t}$ plus the EV charging and discharging does  not exceed the peak capacity $c^{+}_{z,t}$ of the zone. Finally, constraint~\eqref{con_tgt} ensures that the EV energy for vehicle $i$ at the end of the time period $\mathbb{T}$ exceeds a given target $e^\text{tgt}_{i}$.


\section{Analysis design}
\subsection{Scenarios}

Our example scenarios consider a week of EV operation across three zones. The zones are Central Business District (CBD), Suburb and Rural, and have demand profiles based on zones in Melbourne, Australia. We consider both uni-directional and bi-directional charging. For uni-directional charging we replace \eqref{con_power} with $0 \leq p_{i,t} \leq p^{+}$. For bi-directional charging, we allow discharging energy back to the grid to help meet demand and to earn payments for the EV owner. We assume that the discharging price matches the charging price, since discharging offsets other energy use.

Our research primarily focuses on the impact of zone power constraints. In our examples we have a single power constraint for each of three zones, but our model could be extended to include more complicated hierarchies of  constraints. 
\begin{figure}[ht]
\centering
\includegraphics[width=\columnwidth]{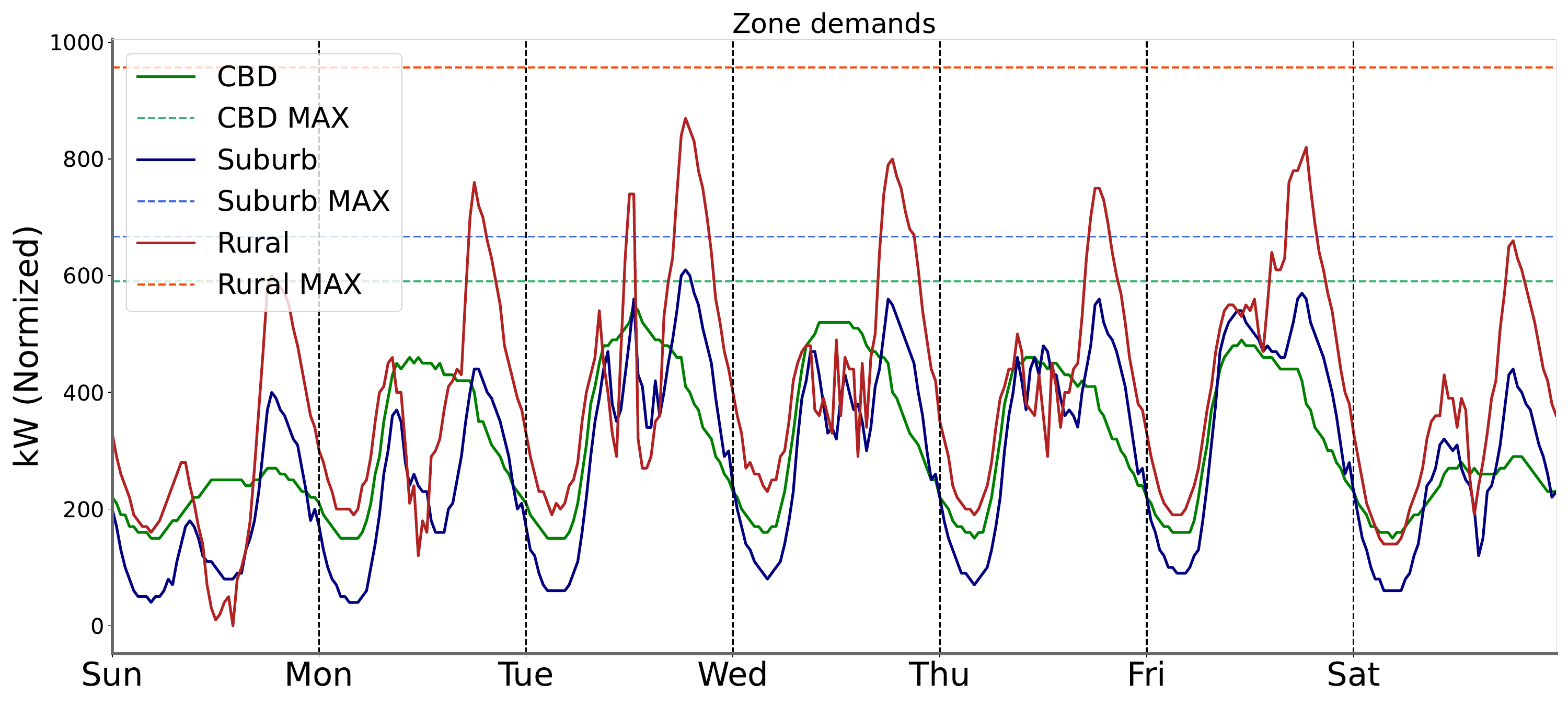}
\caption{Local energy demand for the three example zones}\label{local_dem}
\end{figure}
The local energy demand profiles for each of the three zones were collected from distribution network service providers in Australia~\cite{powercorData} and are shown in Figure~\ref{local_dem}.
We set zone power constraints based on each zone's current power demand. Let  $l^\text{max}_{z}$ be the peak demand of zone $z \in \mathbb{Z}$. The zone power constraints are set to $c^{+}_{z,t} = (1+\eta)l^\text{max}_{z}$ where $\eta$ represents the proportional increase in peak demand. For example, the horizontal dashed lines in Figure \ref{local_dem} shows a 10\% increase in peak demand. In another scenario, we consider the case where only one high-peak zone is constrained. Zones with a higher number of EVs are expected to experience higher energy demand peaks. 

\begin{figure}[ht]
\centering
\includegraphics[width=.93\columnwidth]{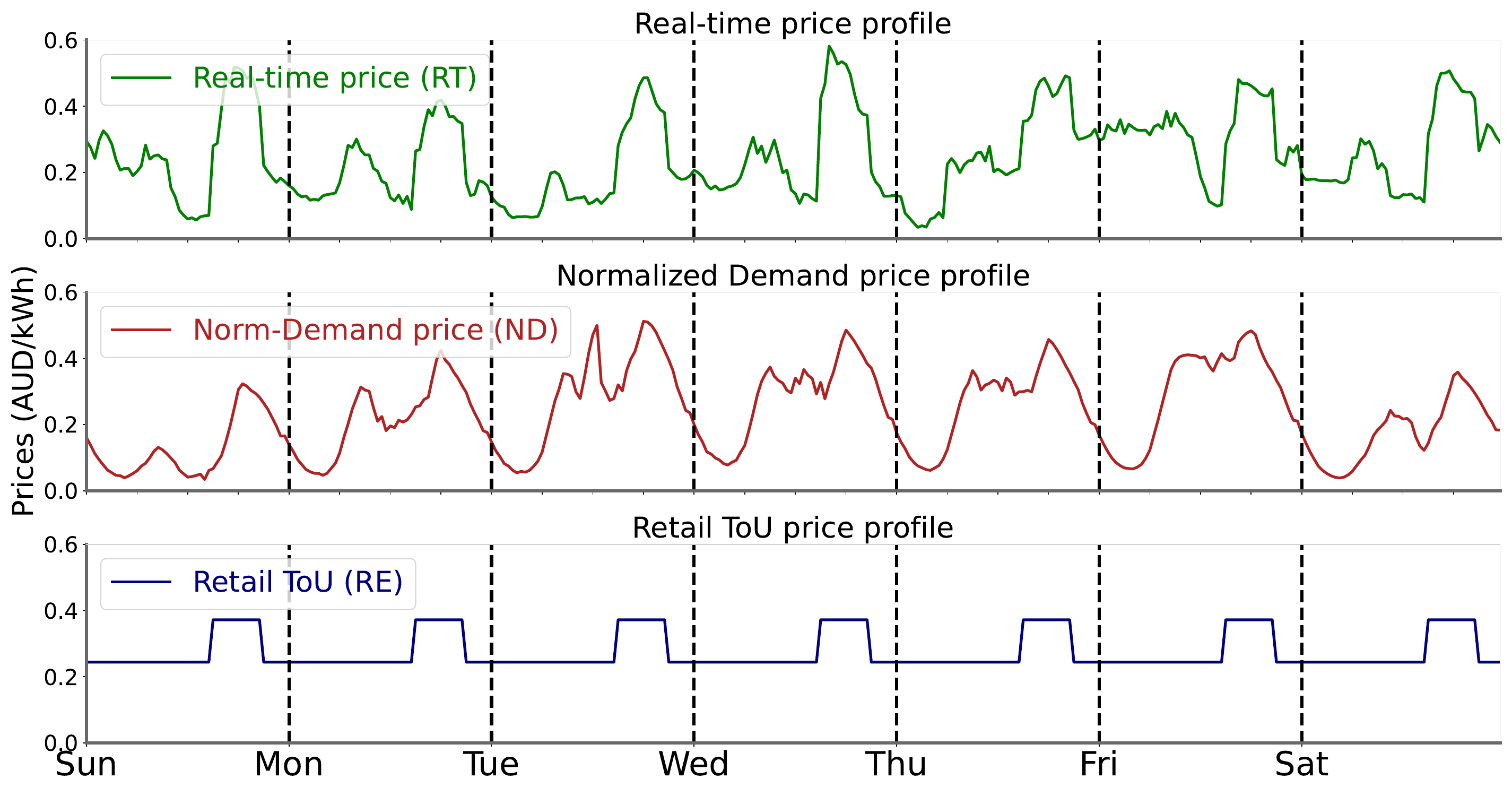}
\caption{Example price profiles}\label{price_dem}
\end{figure}

We also show how different pricing profiles impact EV charging strategies and energy demand in zones. 
We investigate three price profiles, as shown in Figure~\ref{price_dem}:
\begin{itemize}
\item a \textbf{real-time} price profile based on a real-time wholesale price that is the same across all zones, plus a network price that is different for each zone.
\item a \textbf{normalized demand} real-time profile where the price in each zone is proportional to the demand in the zone
\item typical \textbf{retail} ToU tariffs, which are different for each zone.
\end{itemize}

\subsection{Analysis Metrics}

We report four metrics to evaluate the results of each scenario:
\begin{itemize}
\item increase in peak demand in each zone
\item proportion of EV energy use in each zone
\item total charging cost across all EVs
\item proportion of charging energy used for driving.
\end{itemize}

The increase in peak demand in zone $z$ is denoted $\mu_z$, and is given by the ratio of the maximum demand with EVs to the maximum original demand in the zone. The proportion of energy use in zone $z$ is denoted $\xi_z$, and is the ratio of the net EV energy use in zone $z$ to the overall net EV energy use. This metric takes into account the energy used to charge for driving and the energy transferred by EVs between zones.

We evaluate the impact from the EV's perspective by calculating and comparing the total charging cost in each scenario, followed by analyzing the distribution of these costs across EVs to understand the variation. For each EV, the discharged/charged ratio indicates the proportion of the EV energy use that was used for arbitrage. For example, if the ratio is 0.8 then 80\% of the charged energy was discharged back to the grid and 20\% was used for driving and target energy requirement.

\section{Performance Evaluation}
In this section we first describe the experimental setup, focusing on the data used in each of the scenarios. Then, we introduce the experimental results, highlighting the impact of different pricing profiles and the effects of various levels of power capacity constraints on the charging strategies. 

\subsection{Experimental Setup}

Our optimization requires three datasets: local energy demand profiles for each zone, price profiles, and the driving plans of the EVs. 
We consider a week of operation, starting from Sunday, with 48 half-hour intervals. We consider three zones: CBD, Suburb, and Rural zones, based on three network zones in Melbourne, Australia. 
We assume that charging points can handle up to the typical 7.4\,kW in either direction.
Each EV charges from 24\,kWh to a target of 48\,kWh within a battery capacity of 60\,kWh.

Figure \ref{local_dem} shows the normalized local load profiles for the three zones. In high-population zones, local demand is generally expected to be higher. For example, the CBD zone is expected to have greater demand than rural zones. We normalize the local profiles to highlight the distinct characteristics and patterns of each zone's demand profile.

Figure \ref{price_dem} shows the electricity price datasets used in the experiments. The wholesale price profile includes the wholesale energy price from the Australian Energy Market Operator and network costs from various distribution network service providers, while the retail ToU tariffs are sourced from electricity retailers in Victoria, Australia. 

Regarding EV driving plans, we emulate driving behavior across different types of EV users, including ``day-workers'', ``logistics users'', and ``taxi drivers''. We classify destinations for all EVs as ``residential'', ``work'', and ``other,'' and distribute these across the three zones as shown in Table~\ref{gen_ratio}. For example, residential locations for EV users are 10\% in the CBD, 80\% in the Suburb, and 10\% in the Rural zone. 

\begin{table}[t]
\centering
\caption{Distribution of destinations}
\label{gen_ratio}
\begin{tabular}{lccc} 
\toprule
  & CBD & Suburb & Rural  \\ 
\midrule
Residential & 10\% & 80\% & 10\% \\ 
Work        & 70\% & 10\% & 20\% \\ 
Other       & 30\% & 40\% & 30\% \\
\bottomrule
\end{tabular}
\end{table}

\begin{figure}[t]
\centering
\includegraphics[width=.4\columnwidth]{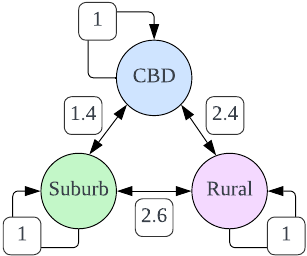}
\caption{The simulated energy consumption (kWh) for traveling within and across zones.} \label{evconsume}
\end{figure}

Each EV's driving plan is paired with the three locations randomly but follows the distribution ratio. However, some EVs may have different destinations within the same zone. To simulate daily driving patterns, EV plans follow a regular schedule based on type, with random settings for parking and commute times. For instance, a day worker drives to the office between 08:00 and 10:00, works 7 to 9 hours, then visits ``other'' location before returning home, with commutes lasting 30 minutes to 1.5 hours. We reference real distances between zones to simulate EV energy consumption for driving. Figure \ref{evconsume} shows the energy consumption required to travel between zones or within the same zone during each time interval $\Delta{t}$.



\subsection{Experimental Results}

\begin{figure}[h]
\centering
\includegraphics[width=\columnwidth]{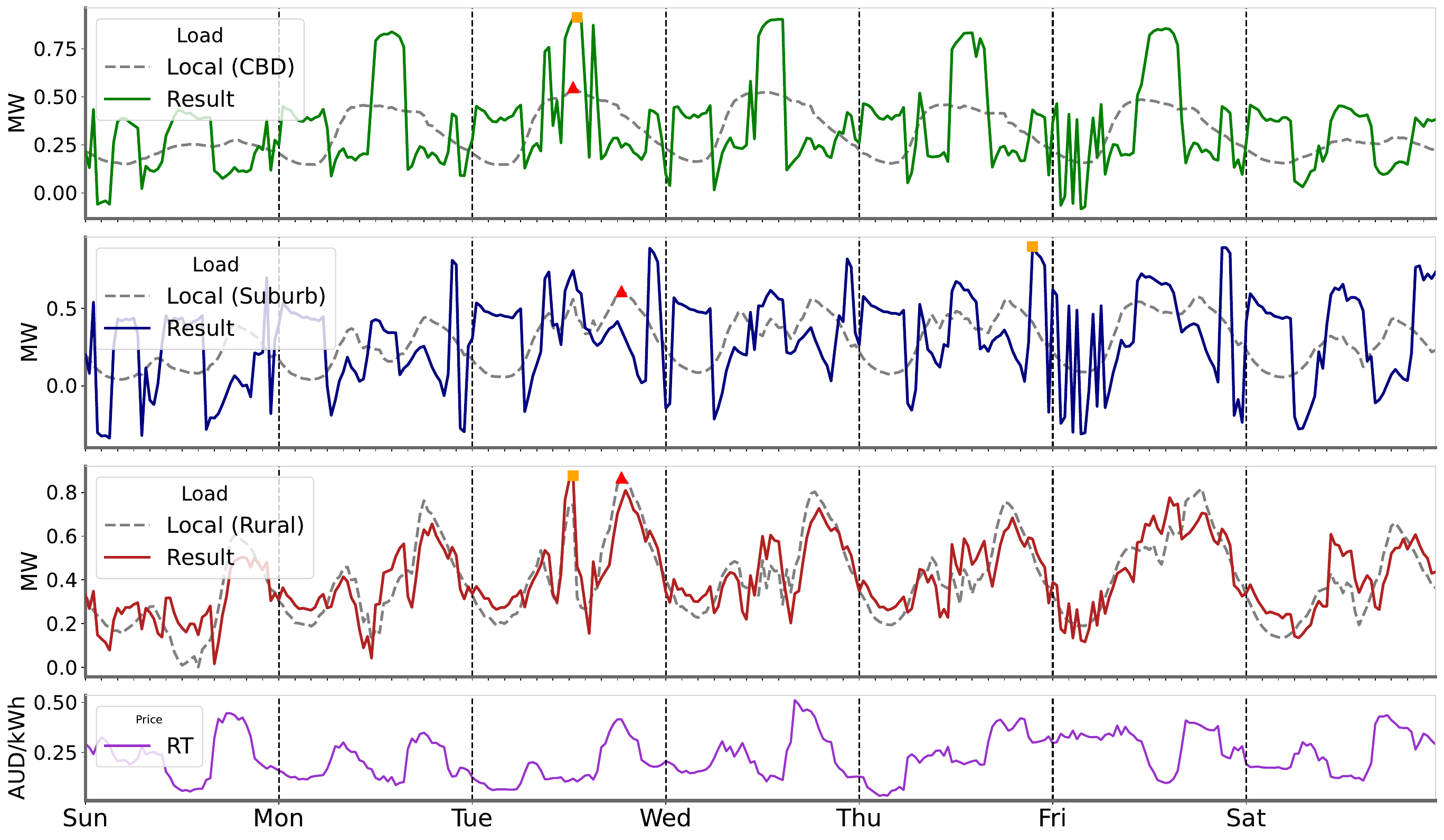}
\caption{Results of the total energy demand versus local energy demand for real-time price profile}\label{default_result}
\end{figure}
\begin{figure}[ht]
\centering
\includegraphics[width=\columnwidth]{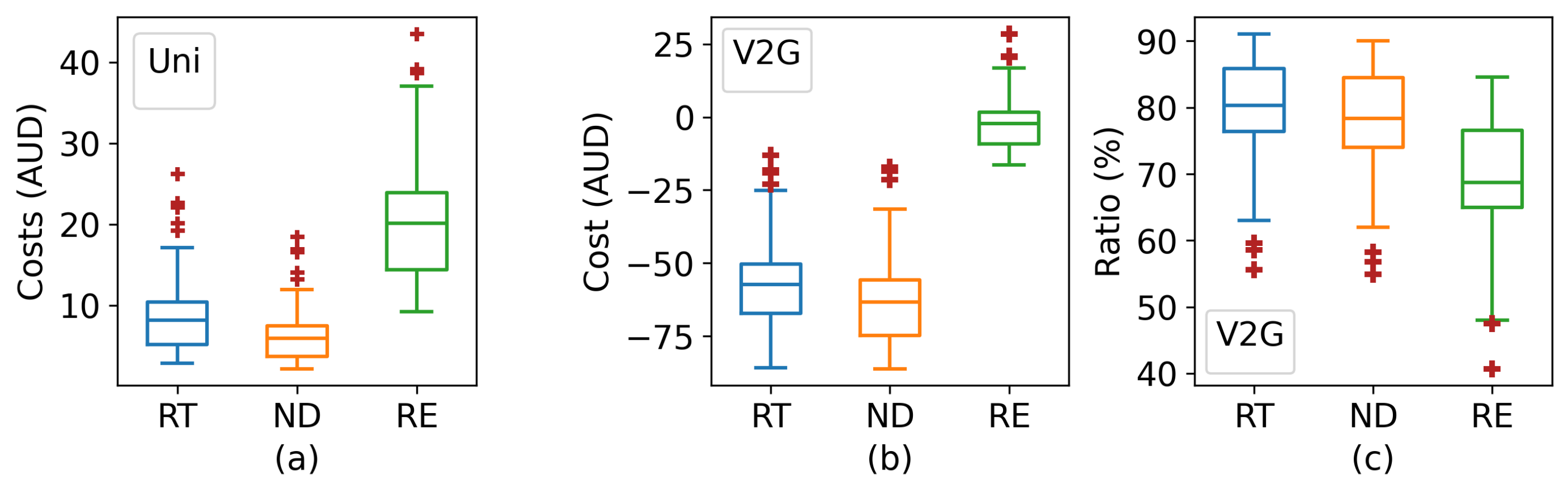}
\caption{Distribution results for different price profiles}\label{box_price}
\end{figure}

Figure~\ref{default_result} illustrates an example of total energy demands from EV charging and local demands across the zones without power constraints. We use Gurobi~11 to solve the Mixed Integer Linear Programming problem. The real-time price profile is applied in this example. Red triangle markers represent the original peak from local demands, while orange circular markers indicate the total energy peak. Table~\ref{tab_price} and Figure~\ref{box_price} display the performance of the metrics in assessing the impact on pricing profiles. In the uni-directional charging scenario, the peak ratio metric indicates an increase in peak demand within the RT and RE profiles, with charged demand primarily concentrated in the CBD and Suburb zones. In the bi-directional (V2G) charging scenario, the results show a notable increase in peak demand compared to the original local peak. In the ND profile, the peak ratios show a less significant increase when pricing is better aligned with demand. This suggests that if prices are more closely matched to demand, the impact of charging demand on the grid would be less substantial. Additionally, the energy ratio in the Suburb zone exhibits higher demand under the ND profile, likely due to its alignment with the energy consumption patterns of primary daytime EV users.
From the EV’s perspective, Figure \ref{box_price}(b) illustrates the distribution of costs, while Figure \ref{box_price}(c) presents the discharged/charged ratio distribution. Given the dynamic nature of the RT and ND profiles, the RE profile results in a higher overall cost distribution. 

\begin{table}[ht]
\centering
\caption{Results for different price profiles}
\label{tab_price}
\begin{tblr}
{
  column{1} = {}{colsep=4pt},
  column{2} = {}{colsep=3pt},
  column{3-9}  = {}{colsep=4pt,r},
  cell{1}{3} = {c=3}{c},
  cell{1}{6} = {c=3}{c},
  cell{1}{9} = {r=2}{},
  cell{3}{1} = {r=3}{},
  cell{6}{1} = {r=3}{},
  vline{3-4,7} = {1}{},
  vline{3-9} = {1-8}{},
  hline{1,3,9} = {-}{0.1em},
  hline{6} = {-}{},
  hline{2} = {3-8}{},
  hline{4-5,7-8} = {2-9}{},
  stretch=0,
}
    &    & Peak Ratio $\mu_{z}$ (\%) &        &       & Energy Ratio $\xi_{z}$ (\%) &        &            & {Total\\ cost} \\
    &    & CBD        & Suburb & Rural & CBD                   & Suburb & Rural &            \\
Uni- & RT & 153.8      & 107.6  & 100.0   & 33.65               & 48.69  & 17.66      & 850     \\
    & ND & 100.0      & 100.0  & 100.0 & 33.78               & 49.03  & 17.19      & 612     \\
    & RE & 108.8      & 100.0  & 100.0 & 30.19               & 52.01  & 17.80      & 2010     \\
V2G & RT & 166.8      & 147.2  & 100.3 & 50.34               & 36.40  & 13.26      & -5774      \\
    & ND & 164.5      & 120.3  & 92.6  & 32.46               & 63.15  & 4.39       & -6305    \\
    & RE & 166.8      & 126.0  & 100.3 & 45.23               & 33.96  & 20.81      & -214   
\end{tblr}
\begin{threeparttable}
\begin{tablenotes}
    \item[Note:] RT: Real-time price, ND: Normalized-Demand price, RE: Retailer ToU price, Uni-: Uni-directional charging, V2G: Bi-directional charging. Total cost in AUD.
\end{tablenotes}
\end{threeparttable}
\end{table}

\begin{figure}[ht]
\centering
\includegraphics[width=.9\columnwidth]{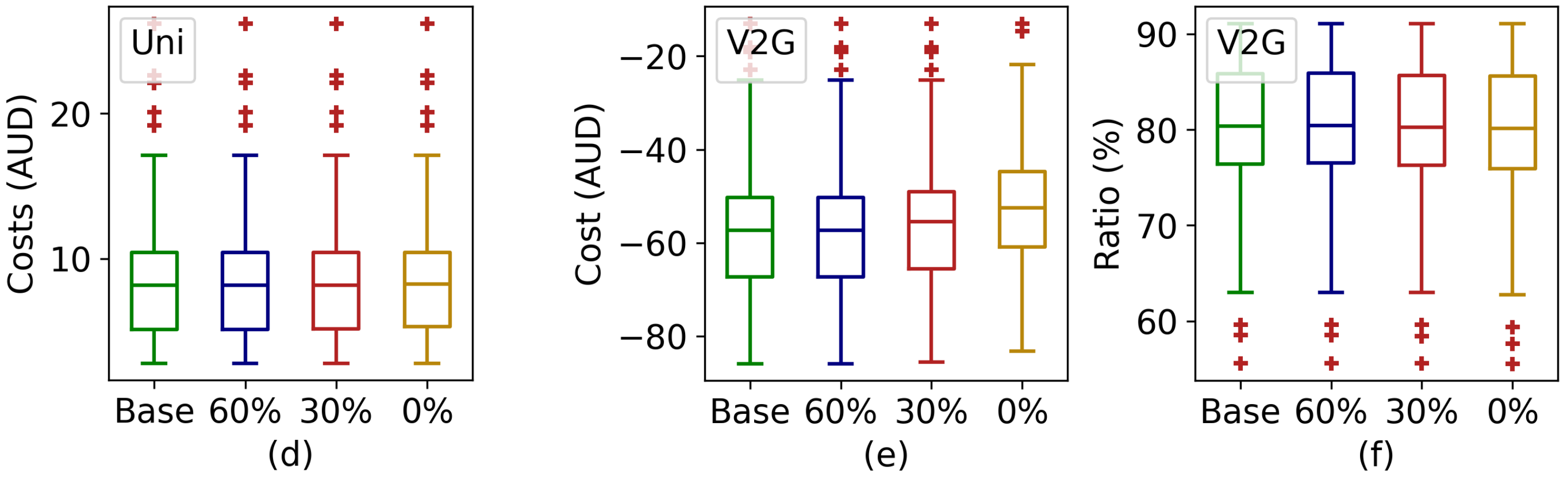}
\caption{Distribution results with zone power constraints (All-zone). (d): Uni-directional, (e): Cost distribution with V2G, (f): Discharged/Charged ratio with V2G.}\label{box_allcon}
\end{figure}

Table \ref{tab_con} and Figure \ref{box_allcon} present the results for the power capacity constraints scenario under the fully constrained condition. In this scenario, we applied the real-time price profile in all experiments. In the uni-directional charging scenario, the peak ratio decreases with varying levels of peak constraint $c^{+}_{z,t}$. Interestingly, despite different constraint levels, the total energy ratio remains consistent. While it was anticipated that stricter constraints would shift charging to other zones, this shift did not occur in the experiments. Additionally, the total EV charging cost and the distribution of EV charging costs in Figure \ref{box_allcon} show no significant increase. A similar trend occurs in the V2G scenario. When the zone power constraint percentage is reduced to 0\%, a clear shift in the demand ratio from the CBD zone to the suburb zone becomes evident, indicating the flexibility of bi-directional charging in reducing charging costs. Moreover, the distribution of EV charging costs shows a noticeably greater increase compared to the uni-directional charging scenario. 

\begin{figure}[h]
\centering
\includegraphics[width=\columnwidth]{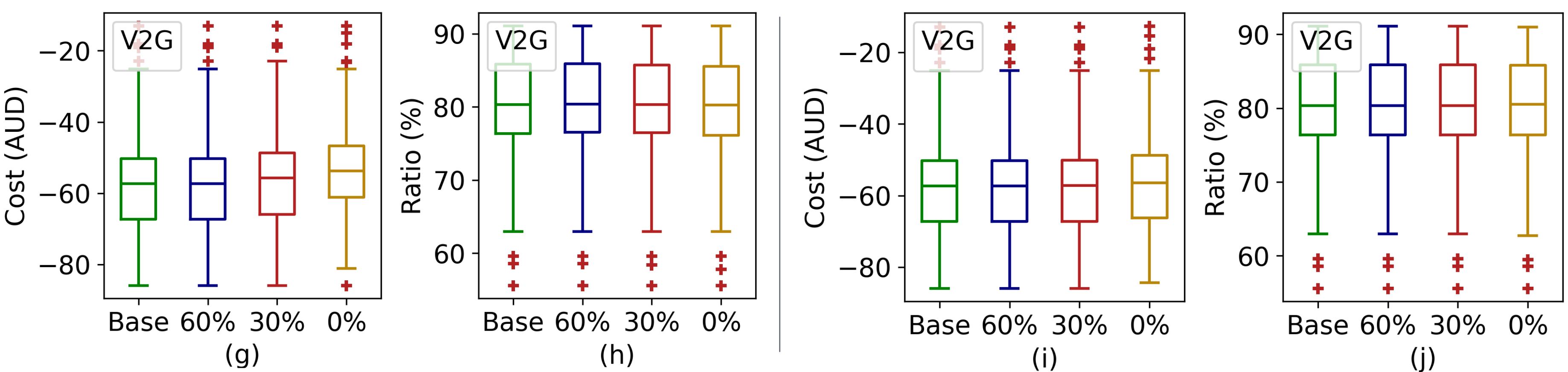}
\caption{Distribution results with zone power constraints. (g),(h): Cost and Discharged/Charged distributions in only-CBD. (i),(j): Cost and Discharged/Charged distributions in only-Suburb.}\label{box_condcon}
\end{figure}

Table~\ref{tab_con} and Figure~\ref{box_condcon} present the results with zone power constraints in the V2G scenario. When the constraint percentage in the CBD decreases, we observe that the peak power in the suburb zone does not increase. A similar outcome occurs in the suburb-constrained scenario. Demand ratios in both scenarios exhibit a comparable pattern: when the power constraint scaled percentage $\eta$ is set to 0, charging demand shifts to other zones. Figure \ref{box_condcon} also indicates a slight increase in charging costs across different constraint levels, while discharged/charged ratios remain largely unaffected by the constraints.

\begin{table}
\centering
\caption{Results with zone power constraints}
\label{tab_con}
\begin{tblr}{
  column{1} = {}{colsep=4pt},
  column{2} = {}{colsep=3pt},
  column{3-9}  = {}{colsep=4pt,r},
  row{7} = {}{rowsep=0.5pt},
  cell{1}{3} = {c=3}{},
  cell{1}{6} = {c=3}{},
  cell{1}{9} = {r=2}{},
  cell{3}{1} = {c=2}{c},
  cell{4}{1} = {r=3}{},
  cell{7}{1} = {c=9}{},
  cell{8}{1} = {c=2}{c},
  cell{9}{1} = {r=3}{},
  cell{12}{1} = {r=3}{},
  cell{15}{1} = {r=3}{},
  vline{4,7} = {1}{},
  vline{3-9} = {1-3,4-6,8-17}{},
  vline{2,4-9} = {3,8}{},
  hline{1,3-4,7-9,12,15,18} = {-}{},
  hline{1,18} = {-}{0.1em},
  hline{2} = {3-8}{},
  hline{5-6,10-11,13-14,16-17} = {2-9}{},
  stretch = 0.1,
}
         &    & Peak Ratio $\mu_{z}$ (\%) &        &       & Energy Ratio $\xi_{z}$ (\%) &        &       & {Total\\cost} \\
         &    & CBD        & Suburb & Rural & CBD                 & Suburb & Rural &            \\
Uni-Base &    & 153.9      & 107.6  & 100.0 & 33.59               & 48.74  & 17.67 & 850     \\
All ($\eta$)      & 60\% & 153.9      & 107.6  & 100.0 & 33.59               & 48.74  & 17.67 & 850     \\
         & 30\% & 130.0        & 107.6  & 100.0 & 33.26               & 49.06  & 17.68 & 851     \\
         & 0\%  & 100.0      & 100.0  & 100.0 & 31.35               & 50.34  & 18.31 & 858     \\
         &    &            &        &       &                     &        &       &            \\
V2G-Base &    & 166.8      & 147.3  & 100.3 & 50.34               & 36.40  & 13.26 & -5774      \\
{All\\($\eta$)}     & 60\% & 160.0        & 147.3  & 100.3 & 50.31               & 36.41  & 13.28 & -5770    \\
         & 30\% & 130.0        & 130.0    & 100.3 & 43.0                  & 43.24  & 13.76 & -5600    \\
         & 0\%  & 100.0      & 100.0  & 100.0 & 30.92               & 51.08  & 18.0    & -5164    \\
{CBD\\($\eta$)}     & 60\% & 160.0        & 147.3  & 100.3 & 50.31               & 36.41  & 13.28 & -5770    \\
         & 30\% & 130.0      & 147.3  & 100.3 & 42.91               & 43.38  & 13.71 & -5608    \\
         & 0\%  & 100.0      & 147.3  & 100.3 & 26.92               & 56.46  & 16.62 & -5282    \\
{Suburb\\($\eta$)}   & 60\% & 166.8      & 147.3  & 100.3 & 50.34               & 36.40  & 13.26 & -5774      \\
         & 30\% & 166.8      & 130.0    & 100.3 & 50.61               & 36.07  & 13.32 & -5767      \\
         & 0\%  & 166.8      & 100.0  & 100.3 & 54.79               & 30.59  & 14.63 & -5658    
\end{tblr}
\begin{threeparttable}
\begin{tablenotes}
    \item[Note:] The Baseline refers to the un-constrained model. Uni-Base: Baseline in uni-directional charging. V2G-Base: Baseline in bi-directional charging. Total cost in AUD.
\end{tablenotes}
\end{threeparttable}
\end{table}

\section{Conclusions and Future Work}
We have developed an optimization model for minimizing the cost of charging a fleet of EVs in scenarios where the EVs move between spatial zones with different power constraints, and illustrated the method for scenarios with different price profiles, zone demand constraints, and with uni-directional and bi-directional charging.

Our example scenarios show that price has an impact on where and when EVs are charged, and that the real-time and normalized-demand price profiles give lower overall cost than the retail pricing for uni-directional charging, and greater overall profit for bi-directional charging. In the zone power constraint scenario, experimental results indicate that as constraints become more restrictive, EV charging demand shifts to less restricted zones. However, this shift results in only a minor increase in overall charging costs.

Minimizing the total cost for a fleet of EVs does not necessarily give a solution that is fair to all EV owners. This could be addressed in further work using fair cost allocation methods from cooperative game theory.




\bibliographystyle{IEEEtran}
\bibliography{11_ref}


\end{document}